\documentclass[11.5pt]{amsart}
\usepackage{graphicx}
\usepackage{amssymb}
\usepackage{amsmath}
\usepackage{amsthm}
\usepackage{amscd}

\newtheorem{thm}{Theorem}[section]
\newtheorem{sub}[thm]{}
\newtheorem{cor}[thm]{Corollary}
\newtheorem{lem}[thm]{Lemma}
\newtheorem{exm}[thm]{Example}

\newtheorem{prop}[thm]{Proposition}
\theoremstyle{definition}

\theoremstyle{remark}

\numberwithin{equation}{section}

\begin{document}

\title[Comodules of $U_q(sl_2)$ and modules of $SL_q(2)$ via quiver ]{Comodules of
$U_q(sl_2)$ and modules of $SL_q(2)$ \\ via quiver}
\author[Xiao-Wu Chen \ and\   Pu Zhang
] {Xiao-Wu Chen$^A$ \ and \  Pu Zhang$^{B,*}$ }

\thanks{$^*$ The corresponding author}
\thanks{Supported in part by the National Natural Science Foundation of China (Grant No. 10271113 and No. 10301033)
and the Doctoral Foundation of the Chinese Education Ministry. The
first named author is supported by the AsiaLink project "Algebras
and Representations in China and Europe" ASI/B7-301/98/679-11}

%
\maketitle

\begin{center}
$^A$Department of Mathematics \\University of Science and
Technology of China \\Hefei 230026, Anhui, P. R. China
\end{center}

\vskip5pt

\begin{center}
$^B$Department of Mathematics\\ Shanghai Jiao Tong University\\
Shanghai 200240, P. R. China
\end{center}
\vskip10pt
\begin{center} xwchen$\symbol{64}$mail.ustc.edu.cn \ \ \ \

pzhang$\symbol{64}$sjtu.edu.cn\end{center}

\begin{abstract}
The aim of this paper is to construct comodules of $U_q(sl_2)$ and
modules of $SL_q(2)$ via quiver, where $q$ is not a root of unity.

By embedding  $U_q(sl_2)$ into the path coalgebra
$k\mathcal{D}^c$, where $\mathcal{D}$ is the Gabriel quiver of
$U_q(sl_2)$ as a coalgebra, we obtain a basis of $U_q(sl_2)$ in
terms of combinations of paths in the quiver $\mathcal{D}$; this
special basis enable us to describe the category of
$U_q(sl_2)$-comodules by certain representations of $\mathcal{D}$;
and this description further permits us to construct a class of
modules of $SL_q(2)$, from certain representations of
$\mathcal{D}$, via the duality between $U_q(sl_2)$ and $SL_q(2)$.
\end{abstract}

\vskip10pt

\section{\bf Introduction}

\vskip10pt

Drinfeld [Dr] has established a duality, between the quantized
enveloping algebra $U_q(sl_2)$ and the quantum deformation
$SL_q(2)$ of the regular function ring on $SL_2$ (see [K], VII).
This has been extended between $U_q(sl_n)$ and $SL_q(n)$ by
Takeuchi [T]. Therefore, any $U_q(sl_n)$-comodule (resp.
$SL_q(n)$-comodule) can be endowed with a $SL_q(n)$-module
structure (resp. a $U_q(sl_n)$-module), in a canonical way (see
e.g. (5.1) below). However, this duality does not give
$U_q(sl_n)$-comodules (resp. $SL_q(n)$-comodules) from
$SL_q(n)$-modules (resp. $U_q(sl_n)$-modules).

Modules of $U_q(\bf\textrm{g})$ have been extensively studied (see
e.g. [L], [Ro], [J]), and it depends on $q$: when $q$ is not a
root of unity, any finite-dimensional module is semi-simple, and
the finite-dimensional simple module is a deformation of a
finite-dimensional simple $\bf\textrm{g}$-module. Another thing of
$U_q(\bf\textrm{g})$ which depends on $q$ is its coradical
filtration ([Bo], [CMus], [M1], [M$\ddot{{\rm u}}$]): when $q$ is
not a root of unity, the graded coalgebra $U_q(\bf\textrm{g})$ is
coradically graded.

The study of $SL_q(n)$-comodules can be also founded, e.g. in
[PW], [CK] (see also [G]). However there are few works on
$U_q(sl_n)$-comodules. A possible reason of this lack might be
that there are no proper tools to construct $U_q(sl_n)$-comodules.
The aim of the present paper is to understand the
$U_q(sl_2)$-comodules by using the quiver techniques.

\vskip10pt

In the representation theory of algebras, quiver is a basic
technique (see [ARS], [Rin]). Recently, it also shows powers in
studying coalgebras and Hopf algebras. For example, one can
construct path coalgebras of quivers, define the Gabriel quiver of
a coalgebra, and embed a pointed coalgebra into the path coalgebra
of its Gabriel quiver (see [CMont], [M3], [CHZ]); after this
embedding one can expect to study the comodules of the coalgebra
by certain locally-nilpotent representations of the quiver (see
[C]); and this makes it possible to see the morphisms, the
extensions, and even the Auslander-Reiten sequences (see e.g.
[Sim]). One can also start from the Hopf quivers of groups to
construct non-commutative, non-cocommutative pointed Hopf algebras
(see [CR]); this makes it possible to classify some Hopf algebras
by quivers, whose bases can be explicitly given (see e.g. [CHYZ],
[OZ]).

\vskip 10pt

Inspired by these ideas,  in this paper, we construct comodules of
$U_q(sl_2)$ and modules of $SL_q(2)$ via quiver, where $q$ is not
a root of unity. By embedding the quantized algebra $U_q(sl_2)$
into the path coalgebra $k\mathcal{D}^c$, where $\mathcal{D}$ is
the Gabriel quiver of $U_q(sl_2)$ as a coalgebra, we obtain a
basis of $U_q(sl_2)$ in terms of combinations of paths in the
quiver $\mathcal{D}$ (Theorem 3.5); this special basis enables us
to describe the category of $U_q(sl_2)$-comodules by certain
locally-nilpotent representations of $\mathcal{D}$ (Theorem 4.3);
in particular, we can list all the indecomposable Schurian
comodules of $U_q(sl_2)$ (Theorem 4.7); and this description
further permits us to construct a class of modules of the quantum
special linear group $SL_q(2)$, from certain locally-nilpotent
representations of $\mathcal{D}$, via the duality between
$U_q(sl_2)$ and $SL_q(2)$ (Theorem 5.2).

\vskip10pt

\section{\bf Preliminaries}

Throughout this paper, let $k$ denote a field of characteristic
zero, and $q$ a non-zero element in $k$ with $q^2 \neq 1$. For a
$k$-space $V$, let $V^*$ denote the dual space. Denote by $\Bbb Z$
and $\Bbb N_0$ the sets of integers and  of non-negative integers,
respectively.

\sub{}\rm  By definition $U_q(sl_2)$ is
an associative $k$-algebra generated by  $E,F, K, K^{-1}$, with
relations (see e.g. [K], p.122, or [J], p. 9)
\begin{align*}
&K K^{-1}=K^{-1}K=1, \\
&K EK^{-1}=q^2E, \quad KFK^{-1}=q^{-2}F,\\
&[E,F]= \frac{K-K^{-1}}{q-q^{-1}}.
\end{align*}
Then $U_q(sl_2)$ has a Hopf structure with (see e.g. [K], p.140)

 \begin{align*}
 \Delta(E)=1 \otimes E + E \otimes K,& \quad \Delta(F)=K^{-1} \otimes F + F \otimes 1 ,\\
\Delta(K)=K \otimes K,  & \quad \Delta(K^{-1})=K^{-1} \otimes K^{-1}, \\
\varepsilon(E)=\varepsilon(F)=0, & \quad
\varepsilon(K)=\varepsilon(K^{-1})=1,\\
 S(K)=K^{-1}, \quad S(K^{-1})=K, & \quad  S(E)=-EK^{-1}, \
S(F)=-KF.
\end{align*}
Note that $U_q(sl_2)$ is a Noetherian algebra without zero
divisors, and it has a basis $\{K^lE^iF^j\ |\  i,j \geq 0, l\in
\Bbb Z\}$ (see e.g. [K], p.123).
\par \vskip 5pt

By definition $SL_q(2)$ is an associative $k$-algebra generated by
$a, b, c, d$, with relations (see e.g. [K], p.84)

\begin{align*}
 ba=qab, &\quad db=qbd,\quad ca=qac, \quad dc=qcd,\quad bc=cb, \\
& ad-da=(q^{-1}-q)bc,\quad  da-qbc=1.
\end{align*}
Then $SL_q(2)$ has a Hopf structure with (see e.g. [K], p.84)

\begin{align*}
 \Delta(a)=a \otimes a +b \otimes c, &\quad \Delta(b)=a \otimes b + b \otimes d, \\
 \Delta(c)=c \otimes a + d \otimes c, &\quad \Delta(d)=c \otimes b + d \otimes d, \\
 \varepsilon(a)=\varepsilon (d)=1,& \quad
 \varepsilon(b)=\varepsilon(c)=0,\\
S(a)= d, \quad S(b)=-qb,& \quad  S(c)=-q^{-1}c, \ \ S(d)=a.
\end{align*}

\vskip 10pt

\sub{}\rm By definition a duality between two Hopf algebras $U$
and $H$ is an algebra map $\psi: H\longrightarrow U^*$, such that
$\phi: U\longrightarrow H^*$ is also an algebra map and has the
property

$$\psi(x)(S_U(u))= \phi(u)(S_H(x))$$
for all $u\in U, x\in H$, where $\phi$ is defined by

$$\phi(u)(x) = \psi(x)(u),$$ and $S_U$ and $S_H$ are respectively the antipodes of $U$ and $H$.
\vskip 5pt Suppose that there exists a duality  between $U$ and
$H$. Then there also exists a duality between $H$ and $U$; and
each $U$-comodule can be endowed with an $H$-module structure, and
also each $H$-comodule can be endowed with a $U$-module. \vskip
5pt

We have the following well-known duality between $U_q(sl_2)$ and
$SL_q(2)$. See Theorem VII.4.4 in [K].

\vskip 5pt

\begin{lem} There is a unique algebra map $\psi: SL_q(2) \longrightarrow U_q(sl_2)^*$ such that

\begin{align*}
&\psi(a)(K^lE^iF^j)=\delta_{i,0}\delta_{j,0} q^l+ \delta_{i,1}
\delta_{j,1}q^l,
 &\quad  \psi(b)(K^lE^iF^j)=\delta_{i,1}\delta_{j,0} q^l,\\
 &\psi(c)(K^lE^iF^j)=\delta_{i,0}\delta_{j,1}q^{-l},
 &\quad \psi(d)(K^lE^iF^j)=\delta_{i,0}\delta_{j,0}q^{-l},
 \end{align*}
 where $\delta_{i, j}$ is the Kronecker symbol. This $\psi$ is
 a duality between $U_q(sl_2)$ and
$SL_q(2)$.
\end{lem}

Note that such a $\psi$ is not injective. This duality was
essentially introduced in [Dr], and has been extended to be a
duality between  $U_q(sl_n)$ and $SL_q(n)$ in [T].

\vskip 10pt

\sub{}\rm\ A quiver $Q=(Q_0, Q_1, s,t)$ is a datum, where $Q$ is
an oriented graph with $Q_0$ the set of vertices and $Q_1$ the set
of arrows, $s$ and $t$ are two maps from $Q_1$ to $Q_0$, such that
$s(a)$ and $t(a)$ are respectively the starting vertex and
terminating vertex of $a\in Q_1$. A path $p$ of length $l$ in $Q$
is a sequence $p= a_l\cdots a_2a_1$ of arrows $a_i, \ 1\le i\le
l,$ such that $t(a_i)=s(a_{i+1})$ for $1 \leq i \leq l-1$. A
vertex is regarded as a path of length $0$. Denote by $s(p)$ and
$t(p)$ the starting vertex and terminating vertex of $p$,
respectively. Then $s(p)= s(a_1)$ and $t(p) = t(a_l)$. If both
$Q_0$ and $Q_1$ are finite sets, then $Q$ is called a finite
quiver. We will not restrict ourselves to finite quivers, but we
assume the quivers considered are countable (i.e., both $Q_0$ and
$Q_1$ are countable sets). For quiver method to representations of
algebras we refer to [ARS] and [Rin].

\vskip10pt

Given a quiver $Q$, define the path coalgebra $kQ^c$ (see [CMon])
as follows: the underlying space has a basis the set of all paths
in $Q$, and the coalgebra structure is given by
$$\Delta(p)=\sum\limits_{\beta\alpha =p}\beta \otimes \alpha$$ and
$$\varepsilon(p)=0\ \ \mbox{if}\ \ l \geq 1, \ \ \mbox{and} \ \ \varepsilon(p)=1\ \  \mbox{if}\ \ l=0$$
for each path $p$ of length $l$.

\vskip10pt

\sub{}\rm \ By a graded coalgebra we mean a coalgebra $C$ with
decomposition $C=\bigoplus\limits_{n\geq 0} C(n)$ of $k$-spaces
such that
$$\Delta(C(n)) \subseteq \sum_{i+j=n} C(i)\otimes C(j), \
\ \ \
 \varepsilon(C(n))=0, \ \ \ \ \forall \ n\ge 1.$$

Let $C$ be a coalgebra. Following [Sw], the wedge of two subspaces
$V$ and $W$ of $C$ is defined to be the subspace

$$V \wedge W: = \{c\in C\ | \ \Delta(c)\in V \otimes C + C \otimes W\}.$$
Let $C_0$ be the coradical of $C$, i.e., $C_0$ is the sum of all
simple subcoalgebras of $C$. Define $C_n: = C_0\wedge C_{n-1}$ for
$n\ge 1$. Then $\{C_n\}_{n\ge 0}$ is called the coradical
filtration of $C$.

\par\vskip 5pt

Recall that a graded coalgebra $C=\bigoplus\limits_{n\geq 0} C(n)$
is said to be coradically graded, provided that
$\{C_n:=\bigoplus\limits_{i\leq n} C(i)\}_{n\ge 0}$ is exactly the
coradical filtration of $C$. It was proved in [CMus], 2.2,  that a
graded coalgebra $C=\bigoplus\limits_{n\geq 0} C(n)$ is
coradically graded if and only if $C_0 = C(0)$ and $C_1 =
C(0)\oplus C(1)$.

\vskip10pt

\sub{}\rm \ Let $M$ be a $C$-$C$-bicomodule over a coalgebra $C$.
Denote by ${\rm Cot}_{C}(M)$ the corresponding cotensor coalgebra
(see [D] for the definition and basic properties). This is a
graded coalgebra with $0$-th component $C$ and $1$-th component
$M$. By Proposition 11.1.1 in [Sw], the coradical of ${\rm
Cot}_{C}(M)$ is contained in $C$. It follows that ${\rm
Cot}_{C}(M)$ is coradically graded if and only if $C$ is
cosemisimple.

\vskip5pt

 Note that a path coalgebra $kQ^c$ is graded with the
length grading, and it is coradically graded, and $kQ^c\simeq {\rm
Cot}_{kQ_0}(kQ_1)$ (see [CMon], or [CR]).

\par\vskip 5pt

We need the following observation.

\begin{prop}\ \ Let $C=\bigoplus\limits_{n\geq 0} C(n)$ be a graded coalgebra. Then

\vskip10pt

(i)\ \  There is a unique graded coalgebra map $\theta: C
\longrightarrow {\rm Cot}_{C(0)}(C(1))$
such that $\theta|_{C(i)}=Id$ for $i=0,1$.\\

(ii)\ \ $\theta (x)= \pi^{\otimes {n+1}}  \circ \Delta^{n} (x)$
for all $x \in C(n+1)$ and $n \geq 1$, where $\pi: C
\longrightarrow C(1)$ is the projection, and $\Delta^{n}=(Id
\otimes \Delta^{n-1}) \circ \Delta$ for all $n \geq 1$, with
$\Delta^{0}= Id$.\\

(iii)\ \ If $C$ is coradically graded, then $\theta$ is
injective.\\

(iv)\ \ If $C(0)$ is cosemisimple, and $\theta$ is injective, then
$C$ is coradically graded.
\end{prop}

\noindent {\bf Proof} \quad Clearly, $C(0)$ is a subcoalgebra and
$C(1)$ is naturally a $C(0)$-$C(0)$-bicomodule, and hence we have
the corresponding cotensor coalgebra ${\rm Cot}_{C(0)}(C(1))$. The
statements (i) and (ii) follow from the universal property of a
cotensor coalgebra (see e.g. [Rad], or [CR]).

For the statement (iii), if $C$ is coradically graded, then $C_1 =
C(0)\oplus C(1)$. It follows that $\theta|_{C_1}$ is injective,
and hence $\theta$ is injective, by a theorem due to Heynemann and
Radford (see e.g. [M2], 5.3.1).

If $C(0)$ is cosemisimple, then ${\rm Cot}_{C(0)}(C(1))$ is
coradically graded. The injectivity of $\theta$ implies that $C$
is a graded subcoalgebra of ${\rm Cot}_{C(0)}(C(1))$. Thus $C$ is
also coradically graded. \hfill $\blacksquare$

\vskip 10pt

\sub{}\rm \ Consider a special case of Proposition 2.7 where
$C(0)$ is a group-like coalgebra (i.e., it  has a basis consisting
of group-like elements; or equivalently, $C(0)$ is cosemisimple
and  pointed). In this case we have $C(0)=kG(C)$, where
$$G(C) : =\{g \in C \ |\ \Delta(g)=g \otimes g,
\varepsilon(g)=1\}.$$  Since $C(1)$ is a $C(0)$-$C(0)$-bicomodule,
it follows that

\begin{align*}
 C(1)=\bigoplus\limits_{g, h \in G}{^h {C(1)}^g},
 \end{align*}
where ${^h {C(1)}^g}=\{c \in C(1)\ | \ \Delta(c)= c \otimes g + h
\otimes c \}$. Define a quiver $Q = Q(C)$ as follows: the set of
vertices is $G$, and there are exactly $t_{gh}$ arrows from vertex
$g$ to vertex $h$, where $t_{gh}=\operatorname{dim}_k {^h
{C(1)}^g}$. Then by the universal property of a cotensor coalgebra
(and hence of a path coalgebra), there is a coalgebra isomorphism
${\rm Cot}_{C(0)} (C(1))\simeq kQ^c$, by identifying the elements
of $G(C)$ with the vertices of $Q$ and a basis of ${^h {C(1)}^g}$
with the arrows from $g$ to $h$.

\vskip5pt

Note that the quiver $Q(C)$ is in general not the Gabriel quiver
of $C$. If the graded coalgebra $C=\oplus_{n\geq 0} C(n)$ is
coradically graded, then $Q(C)$ is exactly the Gabriel quiver of
$C$. For the equivalent definitions of the Gabriel quiver of a
coalgebra we refer to [CHZ], Section 2 (see also [CMon], [M3], and
[Sim]). By Proposition 2.7 we have

\begin{cor}\ \ Assume that $C=\bigoplus\limits_{n \geq 0}C(n)$ is a graded coalgebra
with $C(0)$ group-like. Let $Q(C)$ be the quiver associated to $C$
defined as above. Then

\vskip10pt

(i) \ \ There is a  graded coalgebra map $\theta: \ C
\longrightarrow kQ(C)^c$. \\

(ii) \ \ $\theta$ is injective if and only if $C$ is coradically
graded. In this case, $Q(C)$ is exactly the Gabriel quiver of $C$.
\end{cor}

\vskip10pt

\section{\bf $U_q(sl_2)$ as a subcoalgebra of a path coalgebra}

\vskip10pt

In this section, we embed $U_q(sl_2)$ into the path coalgebra of
the Gabriel quiver $\mathcal{D}$ of $U_q(sl_2)$, and then give a
set of basis of $U_q(sl_2)$ in terms of combinations of paths in
$\mathcal{D}$, where $q$ is not a root of unity.

\vskip5pt

Although bases of $U_q(sl_2)$ are already available, but this new
set of basis of $U_q(sl_2)$ given here, which is in terms of
combinations of paths in $\mathcal{D}$, will enable us to describe
the category of $U_q(sl_2)$-comodules, in terms of
$k$-representations of the quiver $\mathcal{D}$.

\vskip10pt

\sub{}\rm For each non-negative integer $n$, let $C(n)$ be the
subspace of $U_q(sl_2)$ with basis the set $\{ K^lE^iF^j \ | \ i,
j\in \Bbb N_0, i+j =n, \ l \in\Bbb Z\}$. Then
$$U_q(sl_2)= \bigoplus\limits_{n \geq 0} C(n)$$
is a graded coalgebra (see for example Proposition VII.1.3 in [K])
with

$$G(U_q(sl_2)) = \{ K^l\ |\ l\in \Bbb Z\}, \ \ \mbox{and} \ \
C(0)=\bigoplus _{l \in\Bbb Z} kK^l.$$ We have in $C(1)$

\begin{align*}
\Delta(K^{l-1}E)&=K^{l-1} \otimes K^{l-1}E + K^{l-1}E \otimes K^l,\\
\Delta(K^lF)&=K^{l-1}\otimes K^lF + K^lF \otimes K^l.
\end{align*}
Note that $C(1)$ has a set of basis $\{\ K^lE, \ K^lF\ |\ \ l\in
\Bbb Z\}$;

$$^{K^{l_2}} {C(1)}^{K^{l_1}}= 0 \ \ \mbox{for} \ \ (l_1,
l_2)\ne (l, l-1), \ \ l\in\Bbb Z,$$ and that for each $l\in\Bbb Z$
we have

$$^{K^{l-1}} {C(1)}^{K^{l}}= k K^{l-1}E \oplus kK^lF,  \ \  l\in\Bbb Z.$$
Therefore, the quiver of $U_q(sl_2)$ as defined in 2.8 is of the
form \vskip 10pt\vskip 10pt
\begin{center}

\setlength{\unitlength}{.9cm}
\begin{picture}(10,0.7)
\put(0.6,0){\circle*{0.03}}
 \put(0.8,0){\circle*{0.03}}\put(1.0,0){\circle*{0.03}}
 \put(1.2,0){\circle*{0.03}}
\put(1.4,0){\circle*{0.03}} \put(1.6,0){\circle*{0.03}}
\put(1.8,0){\circle*{0.06}} \put(1.9,0.15){\vector(1,0){1}}
\put(1.9,-0.15){\vector(1,0){1}} \put(3,0){\circle*{0.06}}
\put(3.1,0.15){\vector(1,0){1}} \put(3.1,-0.15){\vector(1,0){1}}
\put(4.3, 0){\circle*{0.06}} \put(4.4,0.15){\vector(1,0){1}}
\put(4.4,-0.15){\vector(1,0){1}} \put(5.5,0){\circle*{0.06}}
\put(5.6, 0.15){\vector(1,0){1}} \put(5.6, -0.15){\vector(1,0){1}}
\put(6.7,0){\circle*{0.06}} \put(6.8, 0.15){\vector(1,0){1}}
\put(6.8, -0.15){\vector(1,0){1}} \put(7.9, 0){\circle*{0.06}}
\put(8.1, 0){\circle*{0.03}}
 \put(8.3, 0){\circle*{0.03}}
\put(8.5, 0){\circle*{0.03}} \put(8.7, 0){\circle*{0.03}}
\put(8.9, 0){\circle*{0.03}}\put(9.1, 0){\circle*{0.03}}
\thinlines
\end{picture}

\end{center}
\vskip 10pt\vskip 10pt

\noindent We will denote this quiver by $\mathcal{D}$ in this
paper.

\vskip 10pt

\sub{}\rm \ We fix some notations. Index the vertices of
$\mathcal{D}$ by integers, i.e., $\mathcal{D}_0=\{e_l \ |\ l \in
\Bbb Z\}$; there are two arrows from $e_l$ to $e_{l-1}$ for each
integer $l$. Put $I=\{1, -1\}$ and let $I^n$ be the Cartesian
product (understand $I^0:=\{0\}$). Define
$\mathcal{I}=\bigcup\limits_{n \geq 0} I^n$. For each $v \in
\mathcal{I}$, define $|v|=n$ if $v\in I^n$. Write $v$ as $v =(v_1,
\cdots, v_n),$ where $v_j =1$ or $-1$ for each $j$. For any
integer $l$ and $v \in \mathcal{I}$, define

\begin{align*}
P_l^{(v)}= a_{|v|} \cdots a_2a_1
\end{align*}
to be the concatenated path in $\mathcal{D}$ starting at $e_l$ of
length $|v|$, where the arrow $a_j$ is the upper arrow if $v_j=1$,
and the lower one if otherwise,  $1 \leq j \leq |v|$.

For example, $P_l^{(0)}$ is understood to be the vertex $e_l$;
$P_l^{(1)}$ (resp. $P_l^{(-1)}$) is the upper (resp. lower) arrows
starting at the vertex $e_l$ in $\mathcal{D}$. Clearly,
$$\{P_l^{(v)}= P_{l-|v|+1}^{(v_{|v|})}\cdots P_{l-1}^{(v_2)}P_l^{(v_1)}\ |\ l \in\Bbb Z,
v\in \mathcal{I}\}$$ is the set of
all paths in $\mathcal{D}$.

\vskip 5pt

As an application of Corollary 2.9 we have

\begin{lem}\ \ There is a unique graded coalgebra map $\theta : U_q(sl_2)
\longrightarrow k\mathcal{D}^c$ such that $\theta(K^l)=e_l$,
$\theta(K^{l-1}E)=P_{l}^{(1)},$ \ and \  $\theta
(K^lF)=P_l^{(-1)}$, for each integer $l$.

\vskip 5pt

Moreover, if $q$ is not a root of unity, then $\theta$ is
injective. In this case, $\mathcal{D}$ is the Gabriel quiver of
the coalgebra $ U_q(sl_2)$.
\end{lem}

\noindent{\bf Proof} \quad The existence of $\theta$ follows
directly from Corollary 2.9, and the uniqueness follows from the
universal property of a path coalgebra. Note that if $q$ is not a
root of unity, then the graded coalgebra $U_q(sl_2)=
\bigoplus\limits_{n \geq 0} C(n)$ is coradically graded (see [M1],
or [M2], Question 5.5.6).  \hfill $\blacksquare$

\vskip 10pt

\sub{}\rm \ \ For $v \in I^n\subset\mathcal{I}$,  put

\begin{align*}
T_v :&=\{t \ |\ 1 \leq t \leq n, \quad v_t =1\}, \ \ \ \
\chi(v):=q^{2\sum\limits_{t \in T_v} t}, \ \ \ \mbox{if} \
n\ge 1, \ T_v\ne\emptyset;  \\
\chi(v):& = 1, \ \ \ \ \ \ \mbox{otherwise}.\end{align*}

\vskip 10pt

For each $l\in\Bbb Z, \ n\in\Bbb N_0, \ 0\le i\le n$, set

$$b(l, n, i): = \sum _{v\in I^n, \ |T_v|=i}
\chi(v) P_l^{(v)}\in k\mathcal{D}^c.$$ For example, we have

\begin{align*} b(l, 0, 0)& = e_l, \ \ \ \ \ \ \ \ \ \ \ b(l, 1, 0)  = P_l^{(-1)}, \ \ \ \ \ \ \ \ \ \
b(l, 1, 1) = q^{2}P_l^{(1)},\\
b(l, 2, 0)& = P_l^{(-1, -1)}, \ \ b(l, 2, 2) =  q^{6}P_l^{(1, 1)},\\
b(l, 2, 1)& = q^{2}P_l^{(1, -1)}+ q^{4} P_l^{(-1, 1)}.
\end{align*}

\vskip10pt

The main theorem of this section is

\vskip 5pt

\begin{thm} Assume that $q$ is a not a root of unity.
Then as a coalgebra $U_q(sl_2)$ is isomorphic to  the subcoalgebra
of $k\mathcal{D}^c$ with the set of basis

$$\{b(l, n, i) \ |\ 0\leq i \leq n, \ n\in\Bbb N_0, \ \  l \in\Bbb Z\}.$$
\end{thm}

\vskip10pt

For a non-zero element $q$ in $k$,  and non-negative integers
$n\ge m$, the Gaussian binomial coefficient is defined to be

\begin{align*}
{n \choose m}_q =\frac{n!_q}{m!_q (n-m)!_q}
\end{align*}
where $n!_q:=1_q 2_q\cdots n_q$, \ $0!_q:=1$, \ $n_q:=1+q+ \cdots
+ q^{n-1}$.

\vskip10pt

Given a positive integer $n$, and two vectors $s=(s_0, s_1, \dots,
s_{n-1})$,\ $r=(r_0, r_1, \cdots, r_{n-1})\in \Bbb N_0^{n}$ with
the property

$$s_0\ge s_1\ge \cdots\ge s_{n-1}, \ \ \ r_0\ge r_1\ge\cdots\ge r_{n-1},$$ set

\begin{align*}
c(s, r):= {s_0 \choose s_1}_{q^2} \cdots {s_{n-2} \choose
s_{n-1}}_{q^2} {r_0 \choose r_1}_{q^{-2}} \cdots {r_{n-2} \choose
r_{n-1}}_{q^{-2}} q^{2 \sum\limits_{t=1}^{n-1} r_t(s_{t-1}-s_t)}.
\end{align*}

\vskip10pt

\begin{lem}\ \ Put $E': =K^{-1}E\in U_q(sl_2)$.
Then for any non-negative integers $i$ and $j$, with $n:=i+j\ge
1$, we have

\begin{align*}
\Delta^{n-1} (K^lE'^iF^j)&= \sum_{s, r} c(s, r)\quad  K^{l-s_1-r_1} E'^{s_0-s_1}F^{r_0-r_1}\otimes \cdots \\
&\otimes K^{l-s_{n-1}-r_{n-1}}
E'^{s_{n-2}-s_{n-1}}F^{r_{n-2}-r_{n-1}}\otimes
K^lE'^{s_{n-1}}F^{r_{n-1}}
\end{align*} where the sum runs over all the $r$
and $s$ with $s_0=i$ and $r_0=j$.
\end{lem}

\noindent {\bf Proof} \quad  It suffices to prove the formula for
$n\ge 2$. Note that

\begin{align*}
\Delta(E'^i)&=\Delta(E')^i=(K^{-1} \otimes E' + E' \otimes 1)^i\\
             &=\sum_{s_1=0}^i {i \choose s_1}_{q^2} K^{-s_1} E'^{i-s_1} \otimes E'^{s_1}.
 \end{align*}
 So
 \begin{align*}
 \Delta^{n-1}  (E'^i) &=(Id \otimes \Delta^{n-2})
 (\sum_{s_1=0}^i {i \choose s_1}_{q^2} K^{-s_1} E'^{i-s_1} \otimes E'^{s_1})
 \\&= \sum_{s_1=0}^i {i \choose s_1}_{q^2} K^{-s_1} E'^{i-s_1} \otimes \Delta^{n-2}(E'^{s_1}).
 \end{align*}

\vskip10pt

\noindent By induction we have

\begin{align*} \Delta^{n-1} (E'^i)&=
\sum _{0 \leq s_{n-1} \leq s_{n-2}\leq \cdots \leq s_1 \leq i}  {i \choose s_1}_{q^2}
{s_1 \choose s_2}_{q^2} \cdots {s_{n-2} \choose s_{n-1}}_{q^2} \\
&K^{-s_1}E'^{i-s_1} \otimes K^{-s_2}E'^{s_1-s_2} \otimes \cdots
\otimes K^{-s_{n-1}} E'^{s_{n-2}-s_{n-1}} \otimes E'^{s_{n-1}}.
\end{align*}

\vskip10pt

\noindent Similarly, we have

\begin{align*}
\Delta^{n-1} (F^j) &=
\sum _{0 \leq r_{n-1} \leq r_{n-2}\leq \cdots \leq r_1 \leq j}  {j \choose r_1}_{{q^{-2}}}
{r_1 \choose r_2}_{{q^{-2}}} \cdots {r_{n-2} \choose r_{n-1}}_{{q^{-2}}} \\
&K^{-r_1}F^{j-r_1} \otimes K^{-r_2}F^{r_1-r_2} \otimes \cdots
\otimes K^{-r_{n-1}} F^{r_{n-2}-r_{n-1}} \otimes F^{r_{n-1}}.
\end{align*}

\vskip10pt

\noindent Now the formula follows from $\Delta^{n-1}
(K^lE'^iF^j)=\Delta^{n-1} (K^l)\Delta^{n-1} (E'^i)\Delta^{n-1}
(F^j)$ and the identity
$$ \ \ \ \ \ \ \ \ \ \ \ \ \ \ \ \ \ \ \ E'^mK^{-t} =
q^{2mt}K^{-t}E'^m,\ \ \ \ \ \ m, \ t\in \Bbb N_0. \ \ \ \ \ \ \ \
\ \ \ \ \ \ \ \ \ \ \ \ \ \ \ \ \ \ \ \ \ \ \ \ \blacksquare$$

\vskip 10pt

\sub{}\rm {\bf Proof of Theorem 3.5:} \quad Since $q$ is not a
root of unity, it follows from Lemma 3.3 that there is a coalgebra
embedding $\theta: U_q(sl_2) \longrightarrow k\mathcal{D}^c$. Put
$E':=K^{-1}E$. Then $\{\ K^lE'^iF^j \ | \ i, j\in\Bbb N_0, \
l\in\Bbb Z\ \}$ is a basis of $U_q(sl_2)$. Note that

$$\theta(K^l) = e_l, \ \theta(K^lE') = P_{l}^{(1)}, \ \theta(K^lF) =
P_{l}^{(-1)}.$$ Denote by $\pi$ the projection $U_q(sl_2)
\longrightarrow C(1) \simeq k\mathcal{D}_1$. Then

$$\pi(K^{l-1}E)=P_{l}^{(1)}, \ \pi(K^lF)=P_l^{(-1)}, \ \pi(K^lE'^iF^j) = 0 \ \ \ \mbox{for} \ i+j\ge 2.$$
By Proposition 2.7(ii) we have

\begin{align*}
\theta(K^lE'^i F^j)=\pi^{\otimes n} \circ \Delta^{n-1}
(K^lE'^iF^j)
\end{align*}
where $n=i+j$, and both $i$ and $j$ are positive integers. By
Lemma 3.6 and the definition of $\pi$ we have

\begin{align*}
\theta(K^lE'^iF^j)&=\sum_{s, r} c(s, r) \pi(K^{l-s_1-r_1} E'^{s_0-s_1}F^{r_0-r_1})\cdot \cdots \\
&\cdot \pi(K^{l-s_{n-1}-r_{n-1}}
E'^{s_{n-2}-s_{n-1}}F^{r_{n-2}-r_{n-1}})\cdot
\pi(K^lE'^{s_{n-1}}F^{r_{n-1}})
\end{align*}
where the dot means the concatenation of paths, and the sum runs
over all the vectors $s=(s_0, s_1, \cdots, s_{n-1}),\ r=(r_0, r_1,
\cdots, r_{n-1})\in \Bbb N_0^n$, \ with

$$i = s_0\ge s_1\ge\cdots\ge s_{n-1}, \ \ j = r_0\ge
r_1\ge\cdots\ge r_{n-1},$$ such that for each $t, \ 1\le t\le n$,
either

$$s_{t-1}-s_t = 1, \ r_{t-1}-r_t = 0,$$
or
$$s_{t-1}-s_t = 0, \ r_{t-1}-r_t = 1,$$
where $s_n$ and $r_n$ are understood to be zero.

\vskip 5pt

Now, for such a pair $(s, r)$, define $v= (v_1, \cdots, v_n)\in
I^n$ as follows:

$$v_{n-t+1}= 1, \ \ \mbox{if}\ \ s_{t-1}-s_t = 1, \ r_{t-1}-r_t = 0;$$ and
$$v_{n-t+1}= -1, \ \ \mbox{if}\ \ s_{t-1}-s_t = 0, \ r_{t-1}-r_t = 1,$$
for $1\le t\le n.$ Write $(s, r)$ as $(s, r)= (s_v, r_v)$.

\vskip 5pt

Since $(s_{t-1}+r_{t-1})- (s_t+r_t) = 1$ and $s_n+r_n = 0$, it
follows that $s_t+ r_t =n-t$ for $1\le t\le n-1$. Therefore, we
have

\begin{align*}
\theta(K^lE'^iF^j)&= \sum_{s_v, r_v} c(s_v, r_v) P_l^{(v)}\\& =
\sum_{v\in I^n, \ |T_v|=i} c(s_v, r_v) P_l^{(v)}.
\end{align*}

\vskip 5pt

Note that for $s_v = (i=s_0, s_1, \cdots, s_{n-1})$, any number in
the sequence $i-s_1, \ \cdots, \ s_{n-2}-s_{n-1}, \ s_{n-1}$ is
either $1$ or $0$, and that the number of $1$ in the sequence is
exactly $i$. This implies

$${i \choose s_1}_{q^2} \cdots {s_{n-2} \choose s_{n-1}}_{q^2} =
i!_{q^2}.$$

\vskip 5pt

\noindent In order to compute $c(s_v, r_v)$, let $T_v = \{t_1,
\cdots, t_i\}$, with $1\le t_1 < \cdots < t_i\le n$. By an
analysis on the components of

$$r_v = (j=r_0,  \cdots, \ r_{n-t_i}, \ r_{n-t_i+1}, \cdots,\  r_{_{n-t_{(i-1)}}}, \cdots,
\ r_{t_{_1}}, \cdots, \ r_{n-1}),$$ we observe that $r_{n-t_i} =
r_{n-t_i+1}$ since $v_{t_i} = 1$, and $j=r_0, \cdots, r_{n-t_i}$
are pairwise different. It follows that

$$r_{n-t_i} = j-n+t_i.$$
A similar analysis shows that
$$r_{n-t_x} = j-n + t_x+ (i-x), \ \ x = 1, \cdots, i.$$
It follows that

\begin{align*}\sum\limits_{t=1}^{n-1} r_t(s_{t-1}-s_t) &=\sum\limits_{1\le t\le n-1, \ v_{n-t+1}=1}
r_t \\& = r_{n- t_1} + \cdots + r_{n- t_i} \\& = (t_1 + \cdots +
t_i)-\frac{i(i+1)}{2}.
\end{align*}
This shows

\begin{align*}
c(s_v, r_v)=i!_{q^2} j!_{q^{-2}}\ q^{-i(i+1)}\chi(v),
\end{align*} and hence

\begin{align*}\theta(K^lE'^iF^j)&= i!_{q^2} j!_{q^{-2}} q^{-i(i+1)}\sum\limits_{v\in I^n, \
|T_v|=i}\chi(v)P_{l}^{(v)}\\&= i!_{q^2} j!_{q^{-2}}
q^{-i(i+1)}b(l, n, i)\ \ \ \ \ \ \ \ \ \ \ \ \ \ \ \ \ \ \ \ \ \ \
\ \ \ \ \ \ \ \ \ \ \ \ \ \ \ \ \ \ \ (3.1)\end{align*} for $n =
i+j\geq 2$ and any integer $l$. Thus $U_q(sl_2)\simeq
\theta(U_q(sl_2))$ is spanned by

$$\{b(l, n, i) \ |\ 0\leq i \leq n, \ n\in\Bbb N_0, \ \  l \in\Bbb
Z\},$$ while this set is obviously $k$-linearly independent. This
completes the proof. \hfill $\blacksquare$

\vskip 10pt

\section{\bf Comodules of $U_q(sl_2)$}

\vskip 10pt

In this section, by applying Theorem 3.5 we will characterize the
category of the $U_q(sl_2)$-comodules in terms of the
representations of the quiver $\mathcal{D}$ (see Theorem 4.3), and
then list all the indecomposable Schurian $U_q(sl_2)$-comodules
(see Theorem 4.7), where $q$ is not a root of unity.

\vskip10pt

\sub{}\rm Let $Q$ be a quiver (not necessarily finite). By
definition a $k$-representation of $Q$ is a datum $V=(V_e, f_a; \
e\in Q_0, a\in Q_1)$, where $V_e$ is a $k$-space for each $e\in
Q_0$, and $f_a: V_{s(a)} \longrightarrow V_{t(a)}$ is a $k$-linear
map for each $a\in Q_1$. Set $f_p:=f_{a_l} \circ \cdots \circ
f_{a_1}$ for each path $p=a_l \cdots a_1$, where each $a_i$ is an
arrow, $1\le i\le l$. Set $f_e:= Id$ for $e\in Q_0$. Then $f_p$ is
a $k$-linear map from $V_{s(p)}$ to $V_{t(p)}$. A morphism $\phi:
(V_e, f_a; \ e\in Q_0, a\in Q_1)\longrightarrow (W_e, g_a; \ e\in
Q_0, a\in Q_1)$ is a datum $\phi = (\phi_e; \ e\in Q_0)$ such that

$$\phi_{t(a)}f_a = g_a\phi_{s(a)}$$
for each $a\in Q_1$. Denote by $\operatorname{Rep}(k, Q)$ the
category of the $k$-representations of $Q$. We refer the
representations of quivers to [ARS] and [Rin].

\vskip10pt

A representation $V=(V_e, f_a; \ e\in Q_0, a\in Q_1)$ is said to
be locally-nilpotent, provided that for each $e\in Q_0$ and each
$m \in V_e$, there are only finitely many paths $p$ starting at
$e$ such that $f_p(m)\ne 0$.

\par \vskip 5pt

It was observed by Chin and Quinn that there is an equivalence
between the category of the right $kQ^c$-modules and the category
of the locally-nilpotent representations of $Q$ (see [C]). The
functors can be seen from the following.

 \par \vskip 5pt

For a right $kQ^c$-comodule $(M, \rho)$, define for each $e \in
Q_0$

 \begin{align*}
M_{e}:=\{m \in M\ |\ \rho_0 (m)=m \otimes e\}
 \end{align*}
 where $\rho_0=(Id \otimes \pi_0)\rho$, and $\pi_0: kQ^c \longrightarrow kQ_0$ is the projection.
For every path $p$ there is a unique $k$-linear map $f_p: M_{s(p)}
\longrightarrow M_{t(p)}$, such that for each $m\in M_{s(p)}$
there holds

\begin{align*}
 \rho(m)=\sum_{s(p')=s(p)}f_{p'} (m) \otimes p'
 \end{align*}
where $p'$ runs over all the paths with $s(p') = s(p)$. In this
way we obtain a $k$-representation $(M_e, f_a; \ e\in Q_0, a\in
Q_1)$ of $Q$ satisfying $f_p = f_\beta f_\alpha$ for any path
$p=\beta\alpha$. By construction it is clearly a locally-nilpotent
representation. Note that $M$ is a $kQ_0$-comodule with $\rho_0$.
Since $kQ_0$ is group-like, it follows that we have a
$kQ_0$-comodule decomposition

$$M = \bigoplus\limits_{e\in Q_0}M_e.\eqno(4.1)$$

\par \vskip 5pt

Conversely,  given a locally-nilpotent representation $V=(V_e,
f_{a}; \ e \in Q_0, \alpha \in Q_1)$ of $Q$, define

 \begin{align*}
 M:=\bigoplus_{e \in Q_0} V_e
 \end{align*}
 and $\rho: M \longrightarrow M \otimes kQ^c$ by
 \begin{align*}
 \rho(m):=\sum_{s(p)=e} f_p(m) \otimes p
 \end{align*}
for each $m \in V_e$ (where $f_e$ is understood to be
$\operatorname{Id}$ for $e\in Q_0$). Then $\rho$ is well-defined
since $V$ is locally-nilpotent and $(M, \rho)$ is a right
$kQ^c$-comodule. \par

\sub{}\rm  Keep the notations in 3.2. Given a representation $V =
(V_l, V_a; \ e_l\in \mathcal{D}_0, a\in\mathcal{D}_1)$ of the
quiver $\mathcal{D}$, define $f_{l}^{(v)}:=f_{P_l^{(v)}}$, for
each integer $l$ and $v\in \mathcal{I}$. In particular, $f_l^{(0)}
= Id$.

\par \vskip 5pt

With the help of the representations of a quiver and Theorem 3.5,
we can describe the category of the comodules of $U_q(sl_2)$.

\begin{thm} \ \ Assume that $q$ is not a root of unity.
Then there is an equivalence between the category of the right
$U_q(sl_2)$-comodules and the full subcategory of
$\operatorname{Rep}(k, \mathcal{D})$
whose objects $V=(V_l,\ f_a: \  e_l\in \mathcal{D}_0, \ a\in\mathcal{D}_1)$ satisfies the following conditions:\\

(i) \quad $f_{l-1}^{(1)} \circ f_{l}^{(-1)}=q^2f_{l-1}^{(-1)} \circ f_l^{(1)} $ for all $l \in\Bbb Z$. \\

(ii) \quad For any $m \in V_l$, $f_l^{(v)} (m)=0$ for all but
finitely many $v \in \mathcal{I}$.
 \end{thm}
 \vskip 3pt

 \noindent {\bf Proof} \quad By Theorem 3.5, as a coalgebra $U_q(sl_2)$ is isomorphic to the
 subcoalgebra $\mathcal{C}$
 of path coalgebra  $k\mathcal{D}^c$ with the set of basis

$$\{b(l, n, i): = \sum_{v\in I^n,\  |T_v|=i} \chi(v)P_{l}^{(v)}\ |\
 0\leq i \leq n, \ n\in\Bbb N_0, \ l
\in\Bbb Z\}.$$  For a coalgebra $C$, let $\mathcal{M}^{C}$ denote
the category of the right $C$-comodules. So we have the following
embedding of categories
 \begin{align*}
 \mathcal{M}^{U_q(sl_2)} \simeq \mathcal{M}^{\mathcal{C}} \hookrightarrow \mathcal{M}^{k\mathcal{D}^c}
 \hookrightarrow\operatorname{Rep}(k, \mathcal{D}),
 \end{align*}
 where $\mathcal{M}^{\mathcal{C}} \hookrightarrow \mathcal{M}^{k\mathcal{D}^c}$ since $\mathcal{C}$
 is a subcoalgebra of $k\mathcal{D}^c$, and
 $\mathcal{M}^{k\mathcal{D}^c} \hookrightarrow \operatorname{Rep}(k, \mathcal{D})$
 is the embedding described in 4.1. \par \vskip 5pt

Now, the question is reduced to determine all locally-nilpotent
$k$-representations of quiver $\mathcal{D}$ which are right
$\mathcal{C}$-comodules, via the equivalence described in 4.1.

\vskip5pt

It follows from the definition that a representation $V=(V_l, f_a:
e_l\in \mathcal{D}_0,\ a\in \mathcal{D}_1)$ of quiver
$\mathcal{D}$ is locally-nilpotent if and only if the condition
$(ii)$ is satisfied. Assume that such a $V$ is locally-nilpotent,
then $M=\bigoplus\limits_{l \in \Bbb Z} V_l$ becomes a right
$k\mathcal{D}^c$-comodule via

\begin{align*}
\rho(m)=\sum\limits_{v \in \mathcal{I}} f_{l}^{(v)} (m) \otimes
P_l^{(v)} \in M \otimes k\mathcal{D}^c
\end{align*}
for all $m \in V_l, \ l\in\Bbb Z$.

If for an arbitrary fixed $m\in V_l, \ l\in\Bbb Z$, the element
 $\frac{f_{l}^{(v)}(m)}{\chi(v)}$ only depends on $|v|$ and $|T_v|$,
then we can write

\begin{align*}\rho(m)&=\sum\limits_{n\in\Bbb N_0}\sum\limits_{0\le i\le n} \sum\limits_{v\in I^n, \ |T_v| =i}
\frac{f_{l}^{(v)} (m)}{\chi(v)}  \otimes \chi(v)P_l^{(v)}\\& =
\sum\limits_{n\in\Bbb N_0}\sum\limits_{0\le i\le n}
\frac{f_{l}^{(v)} (m)}{\chi(v)} \otimes (\sum\limits_{v\in I^n, \
|T_v| =i} \chi(v)P_l^{(v)})\\&  \in M \otimes
\mathcal{C},\end{align*} and hence $M$ becomes a right
$\mathcal{C}$-comodule. Conversely, if $M$ becomes a right
$\mathcal{C}$-comodule, then we have

\begin{align*}\rho(m)&=\sum\limits_{n\in\Bbb N_0}\sum\limits_{0\le i\le n} \sum\limits_{v\in I^n, \ |T_v| =i}
\frac{f_{l}^{(v)} (m)}{\chi(v)}  \otimes \chi(v)P_l^{(v)}\\& =
\sum\limits_{n\in\Bbb N_0}\sum\limits_{0\le i\le n} m(n, i)\otimes
(\sum\limits_{v\in I^n, \ |T_v| =i} \chi(v)P_l^{(v)})\end{align*}
for some $m(n, i)\in M$. Since

$$\{\chi(v)P_l^{(v)} \ |\ v\in I^n, \ |T_v| =i, \ 0\le i\le n, \ n\in\Bbb N_0, \ l\in\Bbb Z\}$$
is a set of basis of $k\mathcal{D}^c$, it follows that

$$m(n, i) = \frac{f_{l}^{(v)}(m)}
{\chi(v)},$$ which implies that $\frac{f_{l}^{(v)}(m)}{\chi(v)}$
only depends on $|v|$ and $|T_v|$ for an arbitrary fixed $m\in
V_l, \ l\in\Bbb Z$.

\vskip10pt

Now, the condition $(i)$ implies that for an arbitrary fixed $m\in
V_l, \ l\in\Bbb Z$, the element $\frac{f_{l}^{(v)}(m)} {\chi(v)}$
only depends on $|v|$ and $|T_v|$. Conversely, by taking $v=(-1,
1)$ and $v'=(1, -1)$ in $\mathcal{I}$ we obtain

\begin{align*}
\frac {f_{l}^{(-1,1)}(m)} {\chi((-1, 1))}=\frac
{f_{l}^{(1,-1)}(m)} {\chi((1, -1))},
 \end{align*}
which is exactly the condition $(i)$.  This completes the proof.
\hfill $\blacksquare$

 \vskip 10pt

Theorem 4.3 permits us to explicitly construct some
$U_q(sl_2)$-comodules. In the following $q$ is not a root of
unity.

\vskip10pt

 \begin{exm} \ \
 Let $A$ be the quantum plane generated by $X$ and $Y$ subject to the relation $XY=q^2YX$.
 Let $l$ be an integer and $n$ a non-negative integer. Then for any $A$-module
 $U$ one can define a representation $V= V_{(l, n, U)}$ of quiver $\mathcal{D}$ as follows:

 \begin{align*} &V_j:= U, \ \ &\mbox{if}\ \ l \leq j \leq l+n,\\ &V_j: = 0, \ &\mbox{otherise}; \\
 &f_j^{(1)}:=X,&\mbox{if}\ \ l+1 \leq j \leq l+n,\\
 &f_j^{(1)}:=0,&\mbox{otherwise};
  \\ & f_j^{(-1)}:=Y, \ &\mbox{if}\ \ l+1 \leq j \leq l+n,\\ & f_j^{(-1)}:=   0 & \mbox{otherwise}.
 \end{align*}
 where $l$ is any integer and $n \geq 0$. Then by Theorem 4.3, $V$ induces a right $U_q(sl_2)$-comodule.
 \end{exm}

 \vskip10pt

 \begin{exm} \ \ Let $l$ be an integer and $n$ a non-negative integer.

\vskip5pt

(i) \ \ For each $\lambda \in k$, one can define a representation
$V$ of quiver $\mathcal{D}$ as follows:

 \begin{align*}
 &V_j:=k, &\mbox{ if } \ l\leq j \leq l+n,\\ &
 V_j:=0, & \mbox{ otherwise };\\
 &f_{j}^{(1)}:=1, &\mbox{ if  } \ l+1 \leq j \leq l+n,\\ & f_{j}^{(1)}:=0, & \mbox{ otherwise};\\
 &f_{j}^{(-1)}:=\lambda q^{-2(l+n-j)},& \mbox{ if  } \ l+1 \leq j \leq l+n,\\
& f_{j}^{(-1)}:=0, & \mbox{ otherwise}.
 \end{align*}
Then by Theorem 4.3, $V$ induces a right $U_q(sl_2)$-comodule,
which is denoted by $M_{(l, n, \lambda)}.$

\vskip 10pt

(ii)\ \ Consider the representation $V$ of quiver $\mathcal{D}$
defined by:

\begin{align*}
 &V_j:=k, & \mbox{ if }\ l \leq j \leq l+n,\\& V_j:= 0,& \mbox{ otherwise };\\
 &f_{j}^{(1)}:=0, & \ \forall \ j\in\Bbb Z;
\\&\ f_{j}^{(-1)}:=1, & \forall j \in\Bbb Z.\end{align*} Then by Theorem 4.3, $V$ induces a right
$U_q(sl_2)$-comodule, which is denoted by $M_{(l, n, \infty)}$.
\end{exm}

 \vskip10pt

\sub{}\rm \ \ A finite-dimensional right $U_q(sl_2)$-comodule $(M,
\rho)$ is said to be Schurian, if $\operatorname{dim}_k M_j=1$ or
$0$  for each integer $j$, where $M_j:=\{m \in M\ |\ (Id \otimes
\pi_0)\rho(m)= m \otimes e_j\}$ and $\pi_0$ is the projection from
$k\mathcal{D}^c$ to $k\mathcal{D}_0$.

\vskip10pt

\begin{thm}\ \ When the triple $(l, n, \lambda)$ runs over $\Bbb Z\times \Bbb N_0\times  (k\cup \{\infty\})$,
$M_{(l, n, \lambda)}$ gives a complete list of all pairwise
non-isomorphic, indecomposable Schurian right
$U_q(sl_2)$-comodules, where $q$ is not a root of unity.
\end{thm}

\noindent{\bf Proof}\quad Assume that  $M$ is an indecomposable
Schurian right $U_q(sl_2)$-comodule. Set
$\operatorname{Supp}(M):=\{j \in\Bbb Z\ |\ M_j \neq 0\}$. Let $l$
and $l+n$ be the minimal and the maximal elements in
$\operatorname{Supp}(M)$. Then $\operatorname{Supp}(M)\subseteq
\{l, l+1, \cdots,\ l+n\}$. We claim that
$\operatorname{Supp}(M)=\{l, l+1, \cdots,\ l+n\}$.

Otherwise, there exist a $j_0$ such that $l < j_0 < l+n $ and $j_0
\notin \operatorname{Supp}(M)$. Then by $(4.1)$ we have a
$k\mathcal{D}_0$-comodule decomposition

\begin{align*}
M=(\bigoplus\limits_{j < j_0}M_j)\bigoplus (\bigoplus\limits_{j >
j_0} M_j).
\end{align*}
Since $M_{j_0}= 0$, it follows that this is a
$k\mathcal{D}^c$-comodule decomposition,  and hence it is also a
$U_q(sl_2)$-comodule decomposition, which contradicts to the
assumption.

 \par \vskip 5pt

 Note that each $M_j$ is one-dimensional for $l \leq j \leq l+n$. Set
 \begin{align*}
 a_j:=f_j^{(1)} \quad \mbox{and} \quad b_j:=f_j^{(-1)}, \quad l+1 \leq j \leq l+n.
  \end{align*}
Note that for each $j$, we have $a_j \neq 0$ or $b_j \neq 0$
(otherwise, say $a_{j_0}=b_{j_0}=0$, then we again have a
$U_q(sl_2)$-comodule decomposition  $ M=(\oplus_{j <
j_0}M_j)\oplus (\oplus_{j \geq j_0} M_j)$).

By Theorem 4.3 we have $a_jb_{j+1}=q^2b_{j}a_{j+1}$ for all $j$
with $l+1\le j\le l+n-1$. Now, if some $b_{j_0} = 0$, then all
$b_j = 0$ and all $a_j\ne 0$, and hence $M$ is isomorphic to
$M_{(l, n, 0)}$. If some $a_{j_0} = 0$, then all $a_j = 0$ and all
$b_j\ne 0$, and hence $M$ is isomorphic to $M_{(l, n, \infty)}$.
If $a_j\ne 0\ne b_j$ for all $j$, then $M$ is isomorphic to
$M_{(l, n, \lambda)}$ with $\lambda =
\frac{b_{l+1}}{a_{l+1}}q^{2(n-1)}$.

On the other hand, each $M_{(l, n, \lambda)}$ is indecomposable
since its socle is of one dimension, and they are clearly pairwise
non-isomorphic.\hfill $\blacksquare$

\vskip20pt

\section{\bf A class of $SL_q(2)$-modules}

\vskip10pt

Theorem 4.3 characterizes the category of the right
$U_q(sl_2)$-comodules by a full subcategory of the category of the
$k$-representations of $\mathcal{D}$, where $q$ is not a root of
unity, and $\mathcal{D}$ is the Gabriel quiver of $U_q(sl_2)$ as a
coalgebra. This permits us to construct some left
$SL_q(2)$-modules from some representations of quiver
$\mathcal{D}$, via the duality between $U_q(sl_2)$ and $SL_q(2)$.

\vskip10pt

\sub{}\rm \  Recall that the algebra homomorphism $\psi: SL_q(2)
\longrightarrow U_q(sl_2)^*$  in Lemma 2.3 is given by

\begin{align*} &\psi(a)(K^lE'^iF^j)=\delta_{i,0}\delta_{j,0} q^l+
\delta_{i,1} \delta_{j,1}q^{l-1},
 &\quad  \psi(b)(K^lE'^iF^j)=\delta_{i,1}\delta_{j,0} q^{l-1},\\
 &\psi(c)(K^lE'^iF^j)=\delta_{i,0}\delta_{j,1}q^{-l},
 &\quad \psi(d)(K^lE'^iF^j)=\delta_{i,0}\delta_{j,0}q^{-l},
 \end{align*}
where $E' = K^{-1}E$.

\vskip10pt

Let $(M, \rho)$ be a right $U_q(sl_2)$-comodule. Then $M$ becomes
a left $SL_q(2)$-module via

$$x. m := \sum\psi(x)(m_1) m_0,\eqno(5.1)$$

\vskip5pt \noindent for $x \in SL_q(2), \ m \in M$, where
$\rho(m)=\sum m_0 \otimes m_1\in M\otimes U_q(sl_2)$.

\par \vskip 5pt

Let $\mathcal{C}$ be the subcoalgebra of $k\mathcal{D}^c$ with the
set of basis $$\{  b(l, n, i) = \sum\limits_{v\in I^n,\ |T_v|=i}
\chi(v) P_l^{(v)} \ |\  0 \leq i \leq n, \ \ n\in\Bbb N_0, \ \ l
\in\Bbb Z\}.$$ Identifying $U_q(sl_2)$ with $\mathcal{C}$ via
$(3.1)$, we can evaluate $\psi(a)$, $\psi(b)$, $\psi(c)$ and
$\psi(d)$ on this set of basis of $\mathcal{C}$ via Lemma 2.3.
Since

$$b(l, n, i) = \frac{q^{i(i+1)}}{i!_{q^2} j!_{q^{-2}}}
\theta(K^lE'^iF^j)\ \ \mbox{with} \ \ j = n-i,$$ it follows that
the list of the non-zero values is as follows:

\begin{align*}
\psi(a) (b(l, 0, 0))&= \psi(a) (K^l)=q^l,\\
\psi(a)(b(l, 2, 1))& = \psi(a)(q^{2}K^lE'F)=q^{l+1},\\
\psi(b)(b(l, 1, 1))& = \psi(b)(q^{2}K^lE')=q^{l+1},\\
\psi(c)(b(l, 1, 0))& = \psi(c)(K^lF)=q^{-l},\\
\psi(d)(b(l, 0, 0))& = \psi(d)(K^l)=q^{-l}. \\
\end{align*}

\begin{thm} Let $V=(V_l,\ f_a: \  l\in \mathcal{D}_0, \ a\in\mathcal{D}_1)$ be a $k$-representation of
the quiver $\mathcal{D}$ satisfying the following conditions:\\

(i) \quad $f_{l-1}^{(1)} \circ f_{l}^{(-1)}=q^2f_{l-1}^{(-1)} \circ f_l^{(1)}$ for all $l \in\Bbb Z$. \\

(ii) \quad For any $m \in V_l$, $f_l^{(v)} (m)=0$ for all but
finitely many $v \in \mathcal{I}$, where $f^{(v)}_l =
f_{P^{(v)}_{l}}$, $f^{(0)}_l = Id|_{V_l}$.

\vskip5pt

\noindent Then $M = \bigoplus\limits_{l\in\Bbb Z}V_l$ is a left
$SL_q(2)$-module via

\begin{align*}
&a. m = q^l m + q^{l-1}f_{l}^{(1, -1)}(m),
\\&
b.m = q^{l-1}f_{l}^{(1)}(m), \\& c. m = q^{-l}f_{l}^{(-1)}(m),
 \\& d. m := q^{-l}m,
\end{align*}
for each $m\in V_l, \ l\in\Bbb Z,$ where $q$ is not a root of
unity.
\end{thm}

\noindent{\bf Proof}\ \  By Theorem 4.3 $M =
\bigoplus\limits_{l\in\Bbb Z}V_l$ is a right $U_q(sl_2)$-comodule
via

\begin{align*}\rho(m)= \sum\limits_{n\in\Bbb N_0}\sum\limits_{0\le i\le n}
\frac{f_{l}^{(v)} (m)}{\chi(v)} \otimes b(l, n, i)\end{align*}
where $v$ is a fixed element in $I^n$ with $|T_v| = i$, and $m\in
V_l, \ l\in\Bbb Z$. By $(5.1)$, $M$ becomes a left
$SL_q(2)$-module via

\begin{align*}
 x. m := \sum\limits_{n\in\Bbb N_0}\sum\limits_{0\le i\le n}\psi(x)(b(l, n, i))\frac{f_{l}^{(v)}
 (m)}{\chi(v)}.
\end{align*}
It follows that for each $m\in V_l, \ l\in\Bbb Z,$ we have

\begin{align*}
&a. m = \psi(a) (b(l, 0, 0))m + \psi(a)(b(l, 2,
1))\frac{f_{l}^{(1, -1)}(m)}{\chi(1, -1)}\\ &\quad\quad= q^l m +
q^{l-1}f_{l}^{(1, -1)}(m),
\\&
b.m = \psi(b)(b(l, 1, 1))\frac{f_{l}^{(1)}(m)}{\chi(1)} =
q^{l-1}f_{l}^{(1)}(m),
\\& c. m = \psi(c)(b(l, 1, 0))\frac{f_{l}^{(-1)}(m)}{\chi(0)}=
q^{-l}f_{l}^{(-1)}(m),
 \\& d. m := \psi(d)(b(l, 0, 0))m = q^{-l}m.
\end{align*}
\hfill $\blacksquare$

\vskip10pt

Theorem 5.2 permits us to write out explicitly the following
examples of $SL_q(2)$-modules.

\vskip 10pt

\begin{exm} \ \ Let $A$ be the quantum plane generated by $X$ and $Y$ subject to the relation $XY=q^2YX$,
 and $U$ be a left $A$-module, where $q$ is not a root of unity.
 Let $l$ be an integer and $n$ a non-negative integer.  For any element $u\in U$ and $1\le i\le n+1$,
let $U^{n+1}$ denote the direct sum of the copies of $U$, and
$u_i$ denote  the element in $U^{n+1}$ with the $i$-th component
being $u$ and other components being $0$. Then by Theorem 5.2 and
Example
 4.4,
 the copy $U^{n+1}$ becomes a left $SL_q(2)$-module with the
 following actions :

 \begin{align*} &au_i = q^{i+l-1}u_i + q^{i+l-2}YXu_{i-2} , \ \ & 3 \leq i \leq n+1, \ & \ au_i = 0, \ &\mbox{otherise}, \\
 &bu =q^{i+l-2}Xu_{i-1} ,& 2\leq i \leq n+1,\ \ & \ bu_i = 0, \
 &\mbox{otherise},\\
 &cu_i =q^{-(i+l-1)}Yu_{i-1} ,& 2\leq i \leq n+1,\ \ & \ cu_i = 0, \
 &\mbox{otherise},
  \\ & du_i =q^{-(i+l-1)}u_i, & \forall  i\ \ &.
 \end{align*}
 \end{exm}

\vskip10pt

\begin{exm}\ \ Let $V$ be a $k$-space of dimension $n+1$, \ $n\in\Bbb N_0$, with basis $v_0, v_1, \cdots, v_n$.
Let $l$ be an integer, and $q\in k$ be not a root of unity.

\vskip10pt

(i) \ \ Let $\lambda \in k$. Then by Theorem 5.2 and Example
4.5(i), $V$ becomes a left $SL_q(2)$-module via the following
actions, which is denoted again by $M_{(l, n, \lambda)}$

\begin{align*}
&a.v_i=q^{l+i}v_i + \lambda q^{-2n+l+3i-3} v_{i-2}, &\quad 2 \leq
i \leq n, \\ & a. v_i =q^{l+i}v_i,&
\quad i=0,1,\\
&b.v_i=q^{l+i-1}v_{i-1}, &\quad 1 \leq i \leq n,\\& b.v_0=0,\\
&c.v_i=\lambda q^{-2n+i-l}v_{i-1}, & \quad 1 \leq i \leq n, \\ & c.v_0=0,\\
&d.v_i=q^{-(l+i)}v_i,& \quad \forall\  i.
\end{align*}

 \vskip10pt

(ii)\ \ By Theorem 5.2 and Example 4.5(ii), $V$ also becomes a
left $SL_q(2)$-module via the following actions, which is denoted
again by $M_{(l, n, \infty)}$

\begin{align*}
&a.v_i=q^{l+i}v_i, &\forall\ i,
\\& b.v_i=0, &\forall \ i, \\
&c.v_i=q^{-(l+i)} v_{i-1}, &\quad 1 \leq i \leq n,\\& c.v_0=0,\\
&d.v_i=q^{-(l+i)}v_i, &\quad \forall\ i.
\end{align*}

\vskip10pt

Note that $M_{(l, n, \lambda)}$ with $l\in\Bbb Z, \ n\in\Bbb N_0,
\ \lambda\in k\cup\{\infty\}$ are indecomposable, pairwise
non-isomorphic $SL_q(2)$-modules.
\end{exm}

\vskip10pt

\bibliography{}

\end{document}